\documentclass{article}
\usepackage{amsmath}
\usepackage{amsfonts}

\advance \textheight by \topskip
\advance \textheight by 1pt
\makeatother
\def\bea{\begin{eqnarray}}
\def\ena{\end{eqnarray}}
\def\non{\nonumber}

\def\lar{\longrightarrow}

\def\deg{\hbox{deg}}

\def\tr{\hbox{tr}}

\def\ch{\hbox{ch}}

\def\lar{\longrightarrow}
\def\tpsi{\tilde{\psi}}
\def\tH{\tilde{H}}
\def\e{\text{e}}
\def\tb{\tilde{b}}
\def\uq{U_q(\widehat{sl_2})}
\def\uqs1{U_{\sqrt{-1}}(\widehat{sl_2})}
\def\upqs1{U'_{\sqrt{-1}}(\widehat{sl_2})}
\def\vis1{V_{\sqrt{-1}}(\Lambda_i)}
\def\cq{\mathbb{C}(q)}
\newcommand{\qed}{\hbox{\rule{6pt}{6pt}}}
\newtheorem{prop}{Proposition}
\newtheorem{theorem}{Theorem}

\newtheorem{lemma}{Lemma}
\newtheorem{cor}{Corollary}

\title{
The Space of Local Fields as a Module over the Ring of Local Integrals
of Motion
}

\author{
Atsushi Nakayashiki\thanks{
Faculty of Mathematics,
Kyushu University,
Ropponmatsu 4-2-1, Fukuoka 810-8560, Japan, \quad 
e-mail: 6vertex@math.kyushu-u.ac.jp}
}

\date{}
\begin{document}
\maketitle
\begin{abstract}
The chiral space of local fields in Sine-Gordon or the 
SU(2)-invariant Thirring model is studied as a module over 
the commutative algebra ${\cal D}$ of local integrals of motion.
Using the recent construction of form factors by means of quantum affine
algebra at root of unity due to Feigin et al. we construct a 
${\cal D}$-free resolution of the space of local fields. In general
the cohomologies of the de Rham type complex associated 
with the space of local fields are determined and shown to be the irreducible
representations of the symplectic group $Sp(2\infty)$.
Babelon-Bernard-Smirnov's description of the space of local fields
automatically incorporated in this framework.
\end{abstract}

\section{Introduction}
Let $\uqs1$ be the quantum affine algebra at $q=\sqrt{-1}$,
$x^{\pm}_k$, $b_n$ (a part of) Drinfeld generators,
$\vis1$ the level one representation and $\vis1_{2m+i}$ the subspace
of vectors whose $sl_2$ weight is $2m+i$. 

In \cite{FJKMMT} Feigin et al. proved that the chiral subspace of local fields
in the quantum sine-Gordon model at generic coupling constant and 
that in the SU(2) invariant Thirring model (ITM) are isomorphic
to the space
\bea
&&
A_{2m+i}=
\frac{V_{\sqrt{-1}}(\Lambda_i)_{2m+i}}
{x^{-}_0V_{\sqrt{-1}}(\Lambda_i)_{2m+i+2}+(x^{-}_0)^{(2)}
V_{\sqrt{-1}}(\Lambda_i)_{2m+i+4}},
\quad
m\geq 0,
\non
\ena
where $(x^{-}_0)^{(2)}=(x^{-}_0)^2/(q+q^{-1})$ is the divided power 
of $x^{-}_0$.
The local integral of motion with spin $s$, $s$ being odd,
is identified with $b_s$. They become central at $q=\sqrt{-1}$
and $A_{2m+i}$ becomes a module over the commutative ring
\bea
&&
{\cal D}=\mathbb{C}[b_{-1},b_{-3},\ldots].
\non
\ena

To study the structure of $A_{2m+i}$ as a ${\cal D}$-module 
is an interesting problem because it is a common problem for any 
integrable system. In particular it will give a suitable framework
to compare structures of different integrable systems such as
massive integrable field theories, conformal field theories and 
classical integrable systems \cite{BBS,Z}.

The main aim of this article is to determine the ${\cal D}$-module structure of
$A_{2m+i}$ completely.
More precisely we construct a ${\cal D}$-free resolution of $A_{2m+i}$.
Each term of the resolution is described in terms of the Fock space of
free fermions.

The problem to study the space of local fields as a module over the ring
of local integrals of motion was initiated by 
Babelon-Bernard-Smirnov \cite{BBS}. They studied the restricted sine-Gordon
model and derived the $c=1$ character of the Virasoro algebra. 
The important observation in \cite{BBS} is that the fermionization 
simplifies the description of polynomials giving rise to null vectors.
The description of the space of local fields in terms of the representations
of the quantum affine algebra automatically incorporates the fermions 
and related operators in \cite{BBS}. We remark here that it is difficult
to determine the ${\cal D}$-module structure of $A_{2m+i}$ in the previous 
description by minimal form factors \cite{N2,N3,JMT,JMMT}.

The commutative ring ${\cal D}$ and a module over it
determines a de Rham type complex. Let
\bea
&&
\Omega^{\frac{\infty}{2}}=\mathbb{C}dt_{-1}\wedge dt_{-3}\wedge\cdots,
\non
\ena
be the space of highest degree forms
and $\Omega^{\frac{\infty}{2}-p}$ be the vector space generated by 
differential forms which are obtained from $dt_{-1}\wedge 
dt_{-3}\wedge\cdots$ by removing $p$ $dt_i$'s. 
Then the pair of the vector space and a differential
\bea
&&
C^{\frac{\infty}{2}-p}_{2m+i}
=A_{2m+i}\otimes \Omega^{\frac{\infty}{2}-p},
\label{cochain}
\\
&&
d=\sum_{s=1}^\infty b_{-(2s-1)}\otimes dt_{-(2s-1)},
\label{diff}
\ena
defines a complex. The highest cohomology group is isomorphic to
\bea
&&
\frac{A_{2m+i}}{\sum_{s=1}^\infty b_{-(2s-1)}A_{2m+i}}.
\non
\ena
This space is the space of chiral local fields modulo the action of local 
integrals of motion and its basis gives the minimal set of generators of
$A_{2m+i}$ as a ${\cal D}$-module. 
Thus to determine it is in fact a first step to construct a free resolution of 
$A_{2m+i}$.
In turn, using the free resolution of $A_{2m+i}$ it is possible to 
describe all cohomology groups in terms of the Fock space of free fermions.
We show that those cohomology groups become the irreducible
representations of the symplectic group $Sp(2\infty)$, which is obtained
as the inductive limit of $Sp(2n)$'s. 
This result suggests that the cohomology groups of 
(\ref{cochain}), (\ref{diff}) give the universal
structure of the cohomology groups of affine hyperelliptic Jacobian varieties
\cite{N1,NS1}. 
Thus one can think of the result as an example of the comparison of different
integrable systems mentioned above. 
We shall study this subject in a subsequent paper.

The present paper is organized in the following manner.
After introduction the notations concerning the quantum affine algebra 
at $q=\sqrt{-1}$ and its level one representations are explained in section
two. Then the chiral space of local fields are defined using them.
In section 3 fermions are introduced and the integral expressions of
$x^{-}_0$ and $(x^{-}_0)^{(2)}$ are given in terms of them.
The free resolution of $A_{2m}$ are constructed in section 4.
In section five the de Rham type cohomology groups associated with 
$A_{2m}$ are determined. The space $A_{2m+1}$ is
studied in section 6. We prove that $A_{2m+1}$ is a free ${\cal D}$-module.
In section 7 the cohomology groups are shown to
be the irreducible representations of the symplectic group.

\section{The chiral space}
We recall the results of \cite{FJKMMT}. 
Notations mainly follows that paper.
Let $q$ be an indeterminate. The quantum affine algebra $\uq$ is
the $\cq$ Hopf algebra generated by $x^{\pm}_k$ $(k\in \mathbb{Z})$, 
$b_n (n\in \mathbb{Z}\backslash\{0\})$, $t_1^{\pm1}$, $C^{\pm1}$, $D^{\pm}$
with the following defining relations:

\bea
&&
\text{$C^{\pm1}$ is central,}
\quad
CC^{-1}=C^{-1}C=1,
\non
\\
&&
t_1t_1^{-1}=t_1^{-1}t_1=1,
\quad
t_1b_nt_1^{-1}=b_n,
\quad
t_1x^{\pm}_n t_1^{-1}=q^{\pm2}x^{\pm}_n,
\non
\\
&&
DD^{-1}=D^{-1}D=1,
\quad
Dt_1=t_1D,
\non
\\
&&
Db_nD^{-1}=q^nb_n,
\quad
Dx^{\pm}_nD^{-1}=q^nx^{\pm}_n,
\non
\\
&&
[b_m,b_n]=m\frac{[2m]}{[m]^2}\frac{C^m-C^{-m}}{q-q^{-1}}\delta_{m+n,0},
\non
\\
&&
[b_n,x^{\pm}_k]=\pm\frac{[2n]}{[n]}C^{\frac{n\mp|n|}{2}}x^{\pm}_{k+n},
\non
\\
&&
x^{\pm}_{k+1}x^{\pm}_l-q^{\pm2}x^{\pm}_lx^{\pm}_{k+1}
=q^{\pm2}x^{\pm}_kx^{\pm}_{l+1}-x^{\pm}_{l+1}x^{\pm}_k,
\non
\\
&&
[x^{+}_k,x^{-}_l]=
\frac{C^{-l}\varphi^{+}_{k+l}-C^{-k}\varphi^{-}_{k+l}}{q-q^{-1}},
\non
\ena
where 
\bea
&&
\sum_{k\in \mathbb{Z}}\varphi^{\pm}_{\pm k}z^k
=
t_1^{\pm1}
\exp\Big(
\pm\sum_{n=1}^\infty \frac{q^n-q^{-n}}{n}
\Big)
b_{\pm n}z^n.
\non
\ena
Since we do not use the coproduct in this paper we omit the description it.

In order to specialize $q$ to a complex number one has to define the 
so called integral form.
To this end we need some notations. For $x\in \uq$, the
divided power of $x$ is defined by
\bea
&&
x^{(r)}=\frac{x^r}{[r]!},
\quad
[r]=\frac{q^r-q^{-r}}{q-q^{-1}},
\quad
[r]!=[r]\cdots[2][1].
\non
\ena

Let $A=\mathbb{C}[q,q^{-1}]$. Then the integral form $U_A$ is defined to be
the $A$-subalgebra of $\uq$ generated by
$(x^{\pm}_k)^{(r)}$ $(k\in \mathbb{Z}, r\geq 0)$, 
$b_n$ $n\in \mathbb{Z}\backslash\{0\}$, $t_1^{\pm1}$, $C^{\pm1}$, $D^{\pm1}$.

The algebra $\uqs1$ is then defined by
\bea
&&
\uqs1=\Big(\frac{A}{(q-\sqrt{-1})A}\Big)\otimes_A U_A.
\non
\ena
We denote the subalgebra of $U_A$ without $D^{\pm1}$ by $U'_A$ and similarly
for $\upqs1$.

There are two notable properties which distinguish $\upqs1$ 
from the generic $q$ case. 
They are given by the equations
\bea
&&
x^{\pm}_kx^{\pm}_l=-x^{\pm}_lx^{\pm}_k,
\label{skew}
\\
&&
[b_{2n-1}, \upqs1]=0 \quad \forall n\in \mathbb{Z}.
\label{center}
\ena

Let $\alpha_0, \alpha_1$ and $\Lambda_0, \Lambda_1$ be the simple roots and
the fundamental weights of $\widehat{sl_2}$ respectively.
The Frenkel-Jing realization of the level one integrable highest weight 
representation $V(\Lambda_i)$ of $\uq$ determines the level one 
representation of $\uqs1$ \cite{CJ}. 
Define the integral form $V_A(\Lambda_i)$ by the free $A$-module
\bea
&&
V_A(\Lambda_i)=\oplus_{m\in\mathbb{Z}}A[b_n\,|\, n<0]\e^{\Lambda_i+m\alpha_1}.
\non
\ena
The action of $U_A$ on this space is defined in the following way.
For $P\e^\beta\in V_A(\Lambda_i)$ define the actions of $C$, $\partial$ by
\bea
&&
C(P\e^\beta)=q(P\e^\beta),
\non
\\
&&
\partial(P\e^\beta)=<h_1,\beta>(P\e^\beta),
\non
\ena
where $<h_1,\alpha_1>=2$, $<h_1,\Lambda_j>=\delta_{1,j}$.
For $n<0$ $b_n$ acts by the multiplication by itself and for $n>0$ it acts
by taking the commutator $[b_n,\cdot]$.
Define the grading on $V_A(\Lambda_i)$ as
\bea
&&
\deg \,b_n=n,
\quad
\deg\, \e^{\Lambda_i+m\alpha_1}=-m^2-im.
\non
\ena
Let
\bea
&&
x^{\pm}(z)=\sum_{n\in \mathbb{Z}}x^{\pm}_n z^{-n-1}.
\non
\ena
Then the generators of $U_A$ acts on $V_A(\Lambda_i)$ by the following
formulas
\bea
&&
x^{+}(z)
=
\exp\Big(
\sum_{n=1}^\infty \frac{b_{-n}}{n}z^n
\Big)
\exp\Big(
-\sum_{n=1}^\infty \frac{b_{n}}{n}(qz)^{-n}
\Big)
\e^{\alpha_1}z^\partial,
\non
\\
&&
x^{-}(z)
=
\exp\Big(
-\sum_{n=1}^\infty \frac{b_{-n}}{n}(qz)^n
\Big)
\exp\Big(
\sum_{n=1}^\infty \frac{b_{n}}{n}z^{-n}
\Big)
\e^{-\alpha_1}z^{-\partial},
\label{x-}
\\
&&
t_1=q^{\partial},
\quad
Du=q^{\deg\, u}u \quad\text{for a homogeneous $u\in V_A(\Lambda_i)$}.
\non
\ena

We set
\bea
&&
V_A(\Lambda_i)_{2m+i}=A[b_n\,|\, n<0]\e^{\Lambda_i+m\alpha_1}.
\non
\ena
Then the space
\bea
&&
V_{\sqrt{-1}}(\Lambda_i)=(\frac{A}{(q-\sqrt{-1})A})\otimes_A V_A(\Lambda_i),
\non
\ena
becomes a $\uqs1$-module and
\bea
&&
V_{\sqrt{-1}}(\Lambda_i)_{2m+i}=
(\frac{A}{(q-\sqrt{-1})A})\otimes_A V_A(\Lambda_i)_{2m+i},
\non
\ena
becomes its subspace.

Set
\bea
&&
A_{2m+i}=
\frac{V_{\sqrt{-1}}(\Lambda_i)_{2m+i}}
{x^{-}_0V_{\sqrt{-1}}(\Lambda_i)_{2m+i+2}+(x^{-}_0)^{(2)}
V_{\sqrt{-1}}(\Lambda_i)_{2m+i+4}},
\quad
i=0,1,
\quad
m\geq 0.
\non
\ena
Because of (\ref{center}) it becomes a module over the commutative algebra
\bea
&&
{\cal D}=\mathbb{C}[b_{-1},b_{-3},\ldots].
\non
\ena
The chiral
space of local fields 
in SG-model at generic coupling constant and
that in SU(2)-ITM are both isomorphic to $A_{2m+i}$ 
\cite{FJKMMT,JMT,JMMT}.
In SG case the fields are highest weight vectors with the weight $2m+i$
with respect to a certain quantum group $U_p(\widehat{sl_2})$ for some $p$
and in SU(2)-ITM case the fields are highest weight vectors with weight
$2m+i$ respect to $sl_2$.
 Local integrals of motion with spin $2s-1$ are identified 
with $b_{2s-1}$. Since we consider the chiral space,  only the 
local integrals of motion with negative spins act on it.

\section{Fermions}
We consider the Neveu-Schwarz(NS) and Ramond(R) fermions 
$\{\psi_{2n-1},\psi^\ast_{2n-1}\,|\,n\in \mathbb{Z}\}$ and 
$\{\psi_{2n},\psi^\ast_{2n}\,|\, n\in \mathbb{Z}\}$ respectively.
They satisfy the canonical anti-commutation relations
\bea
&&
[\psi_m,\psi_n]_{+}=[\psi_m^\ast,\psi_n^\ast]_{+}=0,
\quad
[\psi_m,\psi_n^\ast]_{+}=\delta_{m,n}.
\non
\ena
In the following, objects with odd indices are for NS-fermions and those with
even indices are for R-fermions.
The vacuum vectors $|m>$ and $<m|$ are introduced by the following relations
\bea
&&
<m|\psi_n=0 \quad \text{for $n\leq m$},
\qquad
<m|\psi_n^\ast=0 \quad\text{for $n> m$},
\non
\\
&&
\psi_n|m>=0 \quad\text{for $n> m$},
\qquad
\psi_n^\ast|m>=0 \quad\text{for $n\leq m$}.
\non
\ena
These vacuums are related by
\bea
&&
\psi_m^\ast|m-2>=|m>,
\quad
<m-2|\psi_m=<m|.
\non
\ena
The Fock spaces $H_m$, $H_m^\ast$
are constructed from $|m>$ and $<m|$ respectively by the same number of 
$\psi_k$ and $\psi^\ast_l$. The pairing between $H_m$ and $H_m^\ast$ are
defined by the condition
\bea
&&
<m|m>=1.
\non
\ena
The fermion operators are introduced as
\bea
&&
\psi(z)=\sum_{n\in \mathbb{Z}}\psi_{2n+1}z^{-2n-1},
\quad
\psi^\ast(z)=\sum_{n\in \mathbb{Z}}\psi_{2n+1}^\ast z^{2n+1}
\quad
\text{for $NS$-fermions},
\non
\\
&&
\psi(z)=\sum_{n\in \mathbb{Z}}\psi_{2n}z^{-2n},
\quad
\psi^\ast(z)=\sum_{n\in \mathbb{Z}}\psi_{2n}^\ast z^{2n}
\quad
\text{for $R$-fermions}.
\non
\ena
Let $\tb_{2n}$ $(n\in\mathbb{Z})$ satisfy the canonical commutation
relations
\bea
&&
[\tb_{2m},\tb_{2n}]=\delta_{m+n,0}.
\label{canonicalboson}
\ena
We set
\bea
&&
H(b)=\sum_{l=1}^\infty\tilde{b}_{-2l}h_{2l},
\quad
h_{2l}=\sum_n \psi_n\psi^\ast_{n-2l},
\non
\ena
where the summation in $n$ of $h_{2l}$ is taken over odd $n\in\mathbb{Z}$ for
NS-fermions and even $n$ for R-fermions.

The boson-fermion correspondence gives the isomorphism of Fock spaces 
\cite{DJKM}
\bea
&&
H_{2m+i}
\simeq 
\mathbb{C}[\tb_{-2},\tb_{-4},\ldots]\text{e}^{\Lambda_i+m\alpha_1},
\non
\\
&&
a|2m+i>\mapsto <2m+i|\e^{H(b)}a|2m+i>.
\non
\ena
By this correspondence the fermion operators are written as \cite{DJKM}
\bea
&&
\psi(z)=\exp(\sum_{n=1}^\infty \tilde{b}_{-2n}z^{2n})
\exp(-\sum_{n=1}^\infty\frac{\tilde{b}_{2n}}{n}z^{-2n})
\text{e}^{-\alpha}z^{-\partial},
\non
\\
&&
\psi^\ast(z)=\exp(-\sum_{n=1}^\infty \tilde{b}_{-2n}z^{2n})
\exp(\sum_{n=1}^\infty\frac{\tilde{b}_{2n}}{n}z^{-2n})
z^{\partial}\text{e}^{\alpha}
\non
\ena
for both NS and R.

Two bosons $b_{2n}$ and $\tb_{2n}$ are related by
$$
\tb_{2n}=
\left\{
\begin{array}{rl}
\frac{-1}{2}b_{2n},&\quad n\geq 1\\
\frac{(-1)^n}{2n}b_{2n},&\quad n\leq -1.
\end{array}\right.
$$

Then we have
\begin{prop}
By the boson-fermion correspondence $x^{-}_0$ and $(x^{-}_0)^{(2)}$
on $V(\Lambda_i)$ are expressed as
\bea
&&
x^{-}_0=\int \frac{dz}{2\pi i}e^{X(z)}\psi(z),
\label{x-0}
\\
&&
(x^{-}_0)^{(2)}=\frac{-i}{2}\int\int_{|z_1|>|z_2|}
\frac{dz_1}{2\pi i}\frac{dz_2}{2\pi i}\tau(\frac{z_2}{z_1})
\text{e}^{X(z_1)}\text{e}^{X(z_2)}\psi(z_1)\psi(z_2),
\label{x-02}
\ena
where
\bea
&&
X(z)=\sum_{n=1}^\infty \frac{\tb_{-(2n-1)}}{2n-1}z^{2n-1},
\non
\\
&&
\tau(z)=\sum_{n=1}^\infty z^{2n-1}-2\sum_{n=1}^\infty z^{2n},
\non
\ena
and $\psi(z)$ is R-fermion for $i=0$ and NS-fermion for $i=1$.
\end{prop}
\vskip3mm

\noindent
{\bf Remark.} The right hand sides of (\ref{x-0}) and (\ref{x-02})
are similar to the operators appeared in \cite{BBS}. The only difference,
besides the overall constant multiples, is 
the coefficients of even and odd powers of $z$ in $\tau(z)$.
\vskip2mm

\noindent
{\it Proof.} The expression for $x^{-}_0$ follows from (\ref{center})
and the definitions of $\psi(z)$, $X(z)$. Let us prove (\ref{x-02}).
For $n\geq 1$, $b_n$ can be written as
\bea
&&
b_n=n(q^n+q^{-n})\frac{\partial}{\partial b_{-n}}.
\non
\ena
We substitute this into (\ref{x-}). In this description of $x^{-}(z)$ we have
\bea
&&
(x^{-}_0)^{(2)}=\frac{1}{2}\lim_{q\lar i}\frac{(x^{-}_0)^2}{q-i}
=\frac{1}{2}\int\int\frac{dz_1}{2\pi i}\frac{dz_2}{2\pi i}
\frac{d}{dq}\Big(x^{-}(z_1;q)x^{-}(z_2;q)\Big)\vert_{q=i},
\non
\ena
where we write the $q$-dependence of $x^{-}(z)$ explicitly and denote 
$i=\sqrt{-1}$..
By calculation
\bea
&&
\Big(\frac{d}{dq}x^{-}(z;q)\Big)\vert_{q=i}
\non
\\
&&
=
i\sum_{n=1}^\infty b_{-n}i^nz^n x^{-}(z_1;i)
+
x^{-}(z;q)\sum_{n=1}^\infty (-1)^{n-1}(2n-1)z^{-(2n-1)}
\frac{\partial}{\partial b_{-(2n-1)}}.
\non
\ena
Then
\bea
&&
\frac{d}{dq}\Big(x^{-}(z_1;q)x^{-}(z_2;q)\Big)\vert_{q=i}
\non
\\
&&
=
i\sum_{n=1}^\infty b_{-n}i^n(z_1^n+z_2^n) x^{-}(z_1;i)x^{-}(z_2;i)
\non
\\
&&
+x^{-}(z_1;i)x^{-}(z_2;i)
\sum_{n=1}^\infty (-1)^{n-1}(2n-1)(z_1^{-(2n-1)}+z_2^{-(2n-1)})
\frac{\partial}{\partial b_{-(2n-1)}}
\non
\\
&&
-i\tau(\frac{z_2}{z_1})x^{-}(z_1;i)x^{-}(z_2;i).
\non
\ena
Here we use (\ref{center}) and the commutation relations of $b_n$.
Due to (\ref{skew}) 
\bea
&&
\int\int\frac{dz_1}{2\pi i}\frac{dz_2}{2\pi i}
(z_1^n+z_2^n)x^{-}(z_1;i)x^{-}(z_2;i)=0
\quad \forall n\in\mathbb{Z}.
\non
\ena
\qed

\section{Resolution of $A_{2m}$}
In terms of the Fock space of fermions $A_{2m}$ is written as
\bea
&&
A_{2m}\simeq
\frac{{\cal D}\otimes H_{2m}}{x^{-}_0({\cal D}\otimes H_{2m+2})
+(x^{-}_0)^{(2)}({\cal D}\otimes H_{2m+4})}.
\non
\ena
Since $(x^{-}_0)^2=0$ and $[x^{-}_0, (x^{-}_0)^{(2)}]=0$ 
it is possible to define a complex
\bea
&&
\cdots
\lar
\frac{{\cal D}\otimes H_{2m+2}}{(x^{-}_0)^{(2)}({\cal D}\otimes H_{2m+6})}
\lar
\frac{{\cal D}\otimes H_{2m}}{(x^{-}_0)^{(2)}({\cal D}\otimes H_{2m+4})}
\lar 
A_{2m} 
\lar 
0,
\label{resol1}
\ena
where the map 
\bea
&&
\frac{{\cal D}\otimes H_{2n}}{(x^{-}_0)^{(2)}({\cal D}\otimes H_{2n+4})}
\lar
\frac{{\cal D}\otimes H_{2n-2}}{(x^{-}_0)^{(2)}({\cal D}\otimes H_{2n+2})}
\non
\ena
is given by the left multiplication by $x^{-}_0$.

\begin{prop}
The module
\bea
&&
\frac{{\cal D}\otimes H_{2n}}{(x^{-}_0)^{(2)}({\cal D}\otimes H_{2n+4})}
\non
\ena
is a free ${\cal D}$-module.
\end{prop}
\vskip2mm

\noindent
{\it Proof.}  
To prove the proposition we write $x^{-}_0$ and $(x^{-}_0)^{(2)}$ 
in the component form.
Let us write
\bea
&&
\e^{X(z)}=\sum_{n=0}^\infty T_{-n}z^n.
\non
\ena
Notice that $T_{-(2n-1)}$ is a homogeneous polynomial of $\tb_{-(2m-1)}$ with
the degree $-(2n-1)$ and it has the form
\bea
&&
T_{-(2n-1)}=\frac{\tb_{-(2n-1)}}{2n-1}+\cdots,
\non
\ena
where $\cdots$ part does not contain $\tb_{-(2n-1)}$. In particular
\bea
&&
{\cal D}=\mathbb{C}[T_{-1},T_{-3},\ldots].
\non
\ena
By calculation we have
\bea
&&
x^{-}_0=\sum_{n=1}^\infty T_{-(2n-1)}\psi_{2n},
\non
\\
&&
(x^{-}_0)^{(2)}=\frac{-i}{2}\sum_{n=1}^\infty
\Big(
\psi_{2(-n+1)}+\sum_{l=1}^\infty Q_{n-1,l}(T)\psi_{2(-n+1+l)}
\Big)\psi_{2n},
\non
\ena
where
\bea
&&
Q_{n,l}(T)=\sum_{n_1+n_2=l, 0\leq n_1, 0\leq n_2\leq n}
(T_{-2n_1}T_{-2n_2}-2T_{-(2n_1+1)}T_{-(2n_2-1)}).
\non
\ena
We omit the tensor symbol $\otimes$ in writing elements of 
${\cal D}\otimes H_{2m}$ etc. for the sake of simplicity.
Using this expression we make a base change of 
${\cal D}\otimes H_{2n}$.
Let us define $\{\tpsi_{2n},\tpsi^\ast_{2n}\}$ in the following way.
First we set
\bea
&&
\tpsi_{2n}=\psi_{2n}\quad \text{for $n\geq 1$},
\non
\\
&&
\tpsi_{-2n}=\psi_{-2n}+\sum_{l=1}^\infty Q_{n,l}(T)\psi_{2(-n+l)}
\quad \text{for $n\geq 0$},
\non
\ena
and write it as
\bea
&&
\tpsi_i=\sum_j B_{ij}\psi_j,
\non
\ena
where $i,j$ are even integers.
Then the matrix $B=(B_{ij})$ is a triangular matrix such that its 
diagonal entries are all $1$ and $B_{ij}=\delta_{ij}$ for $i\geq 2$. 
Thus $B^{-1}$ exists.
Define the matrix $C$ by
\bea
&&
C=(c_{ij})={}^t(B^{-1})
\non
\ena
and set
\bea
&&
\tpsi_i^\ast=\sum_jc_{ij}\psi_j^\ast.
\non
\ena
Then $\{\tpsi_k, \tpsi_l^\ast\}$ satisfy 
the canonical anti-commutation relations
\bea
&&
[\tpsi_k,\tpsi_l]_{+}=[\tpsi_k^\ast,\tpsi_l^\ast]_{+}=0,
\quad
[\tpsi_k,\tpsi_l^\ast]_{+}=\delta_{k,l}.
\non
\ena
Moreover the vacuum is the same for $\{\tpsi_k, \tpsi_l^\ast\}$, that is, 
the following relations hold
\bea
&&
\tpsi_{2n}|2m>=0\quad \text{for $n>m$},
\quad
\tpsi^\ast_{2n}|2m>=0\quad \text{for $n\leq m$},
\non
\\
&&
<2m|\tpsi_{2n}=0\quad \text{for $n\leq m$},
\quad
<2m|\tpsi^\ast_{2n}=0\quad \text{for $n>m$}.
\non
\ena
Therefore we denote the vacuum for $\{\tpsi_k, \tpsi_l^\ast\}$ 
by the same symbol
as for $\{\psi_k, \psi_l^\ast\}$. 
We denote the Fock space of $\{\tpsi_k, \tpsi_l^\ast\}$
by $\tH_{2m}$, $\tH_{2m}^\ast$. Then we have the isomorphism of ${\cal D}$
modules
\bea
&&
{\cal D}\otimes H_{2m}\simeq {\cal D}\otimes \tH_{2m}.
\non
\ena
In terms of $\{\tpsi_k, \tpsi_l^\ast\}$, $x^{-}_0$ and $(x^{-}_0)^{(2)}$ 
take simple forms
\bea
&&
x^{-}_0=\sum_{n=1}^\infty T_{-(2n-1)}\tpsi_{2n},
\non
\\
&&
(x^{-}_0)^{(2)}=\frac{i}{2}\sum_{n=1}^\infty
\tpsi_{-2(n-1)}\tpsi_{2n}.
\non
\ena
In particular $(x^{-}_0)^{(2)}$ defines a map
\bea
&&
(x^{-}_0)^{(2)}:\tH_{2n}\lar \tH_{2n-4}.
\label{x02}
\ena
Thus we have the isomorphism
\bea
&&
\frac{{\cal D}\otimes H_{2n}}{(x^{-}_0)^{(2)}({\cal D}\otimes H_{2n+4})}
\simeq
{\cal D}\otimes \frac{\tH_{2n}}{(x^{-}_0)^{(2)}\tH_{2n+4}}.
\non
\ena
This proves the proposition.
\qed
\vskip2mm

\noindent
{\bf Remark.} The argument used in the proof of the proposition is essentially
due to Babelon-Bernard-Smirnov\cite{BBS}.
\vskip2mm

We set 
\bea
&&
W_{2n}=\frac{\tH_{2n}}{(x^{-}_0)^{(2)}\tH_{2n+4}}.
\non
\ena

\begin{theorem}\label{resol}
The complex (\ref{resol1})
\bea
&&
\cdots
\lar
{\cal D}\otimes W_{2m+2}
\lar
{\cal D}\otimes W_{2m}
\lar 
A_{2m} 
\lar 
0,
\non
\ena
is a ${\cal D}$-free resolution of $A_{2m}$.
\end{theorem}

The following lemmas are sufficient to prove the theorem.

\begin{lemma}\label{lem1}
The complex 
\bea
&&
\cdots
\lar
{\cal D}\otimes \tH_4
\lar
{\cal D}\otimes \tH_2
\lar
{\cal D}\otimes \tH_0
\lar
0
\non
\ena
is exact at $\tH_{2n}$ $n\geq 1$, where the maps are defined by 
the multiplication by $x^{-}_0$.
\end{lemma}

\begin{lemma}\label{lem2}
The map $(\ref{x02})$ is injective for $n\geq 1$.
\end{lemma}

Assuming these lemmas we give a proof of the theorem first.

\noindent
{\it Proof of Theorem \ref{resol}.}
\par

It is sufficient to prove the exactness at
${\cal D}\otimes W_{2n}$, $n>m$. 
Suppose that $v\in \tH_{2n}$ satisfies
\bea
&&
x^{-}_0v=(x^{-}_0)^{(2)}w
\non
\ena 
for some $w\in \tH_{2n+2}$. Then 
\bea
&&
(x^{-}_0)^{(2)}(x^{-}_0w)=0.
\non
\ena
Since $(x^{-}_0)^{(2)}$ is injective by Lemma \ref{lem2},
$x^{-}_0w=0$. By Lemma \ref{lem1}, 
$w=x^{-}_0u$ for some $u\in \tH_{2n+4}$. 
Then
\bea
&&
x^{-}_0(v-(x^{-}_0)^{(2)}u)=0,
\non
\ena
and again by Lemma \ref{lem1}
\bea
&&
v-(x^{-}_0)^{(2)}u=x^{-}_0y
\non
\ena
for some $y\in \tH_{n+2}$. This proves the theorem.
\qed

\vskip2mm
\noindent
{\it Proof of Lemma \ref{lem1}.}
\par
Notice that
\bea
&&
\tH_{2m}=
\sum\mathbb{C}
\tpsi_{2r_1}^\ast\cdots \tpsi_{2r_{k+m}}^\ast
\tpsi_{2s_1}\cdots\tpsi_{2s_k}|0>,
\non
\ena
where the summation is taken for all
\bea
&&
0<r_1<\cdots<r_{k+m},
\quad
s_1<\cdots<s_k\leq 0.
\label{index}
\ena
We set
\bea
&&
d=x^{-}_0=\sum_{s=1}^\infty T_{-(2s-1)}\psi_{2s},
\non
\ena
for the sake of simplicity.
For each $(s_1,\ldots,s_k)$ satisfying (\ref{index}) we set
\bea
&&
\tH_{2m}(s_1,\ldots,s_k)=
\sum_{0<r_1<\cdots<r_{k+m}}
\mathbb{C}
\tpsi_{2r_1}^\ast\cdots \tpsi_{2r_{k+m}}^\ast
\tpsi_{2s_1}\cdots\tpsi_{2s_k}|0>.
\non
\ena
Then $\tH_{2m}$ is a direct sum of $\tH_{2m}(s_1,\ldots,s_k)$'s and
\bea
&&
d\Big(\tH_{2m}(s_1,\ldots,s_k)\Big)\subset \tH_{2m-2}(s_1,\ldots,s_k).
\non
\ena
Thus it is sufficient to prove the exactness of the complex 
\bea
&&
({\cal D}\otimes \tH_{2m}(s_1,\ldots,s_k),d)
\non
\ena
at $m\geq 1$ for each $(s_1,\ldots,s_k)$. 
Set, for $0<r_1<\cdots< r_{k+m}$,
\bea
&&
e_{r_1\ldots r_{k+m}}=
\tpsi_{2r_1}^\ast\cdots \tpsi_{2r_{k+m}}^\ast
\tpsi_{2s_1}\cdots\tpsi_{2s_k}|0>.
\non
\ena
They form a basis of $\tH_{2m}(s_1,\ldots,s_k)$. Now suppose that
\bea
&&
v\in {\cal D}\otimes \tH_{2m}(s_1,\ldots,s_k),
\quad
dv=0,
\non
\ena
and write
\bea
&&
v=\sum_{0<r_1<\cdots< r_{k+m}}
P_{r_1\ldots r_{k+m}}e_{r_1\ldots r_{k+m}},
\quad
P_{r_1\ldots r_{k+m}}\in {\cal D}.
\non
\ena
If we set $t_s=T_{-(2s-1)}$ we have
\bea
&&
d(e_{r_1\ldots r_{k+m}})=
\sum_{i=1}^{k+m}(-1)^{i-1}t_{r_i}
e_{r_1\ldots r_{i-1},r_{i+1}\ldots r_{k+m}}.
\non
\ena
Let $N_1$ be the maximum among $r_1,...,r_{k+m}$ such that 
$P_{r_1\ldots r_{k+m}}\neq0$, $N$ an integer such that $N>N_1$ and
$P_{r_1\ldots r_{k+m}}\in \mathbb{C}[t_1,\ldots,t_{N}]$ for all
$P_{r_1\ldots r_{k+m}}\neq0$. We set
\bea
&&
B=\mathbb{C}[t_1,\ldots,t_N],
\quad
K_p=\sum_{1\leq i_1<\cdots<i_p\leq N}Be_{i_1\ldots i_p}.
\non
\ena
Then 
\bea
&&
v\in K_{m+k}
\non
\ena
and the complex
\bea
&&
0
\lar
K_N
\stackrel{d}{\lar}
K_{N-1}
\stackrel{d}{\lar}
\cdots
\stackrel{d}{\lar}
K_0
\lar
0
\non
\ena
is the Koszul complex associated with the regular sequence $(t_1,\cdots,t_N)$
of $B$. Thus the complex $(K_{\cdot},d)$ is exact at 
$K_n$, $n\geq 1$ \cite{Mat}. 
Since $m\geq 1$, this proves the lemma.
\qed
\vskip2mm

\noindent
{\it Proof of Lemma \ref{lem2}}
\par
Let us set
\bea
&&
\omega=2i (x^{-}_0)^{(2)}=-\sum_{n=1}^\infty \tpsi_{-2(n-1)}\tpsi_{2n}.
\non
\ena
Suppose that $v\in \tH_{2m}$ satisfies 
\bea
&&
\omega v=0.
\non
\ena
Let us write
\bea
&&
v=\sum c(r_1,\ldots,r_k|s_1,\ldots,s_{k+m})
\tpsi_{2r_1}\cdots \tpsi_{2r_k}
\tpsi^\ast_{2s_1}\cdots \tpsi^\ast_{2s_{k+m}}|0>,
\non
\ena
where 
\bea
&&
r_1<\cdots<r_{k}\leq 0<s_1<\cdots<s_{k+m}.
\non
\ena
Let $N$ be the maximum number among $-r_1,\ldots,-r_k,s_1,\ldots,s_{k+m}$ 
such that $c(r_1,\ldots,r_k|s_1,\ldots,s_{k+m})\neq 0$.
Then
\bea
&&
\omega v=\omega_N v,
\quad
\omega_N=-\sum_{n=1}^N \tpsi_{-2(n-1)}\tpsi_{2n}.
\non
\ena
Define
\bea
&&
\eta_N=\sum_{n=1}^N\tpsi_{-2(n-1)}^\ast\tpsi_{2n}^\ast,
\quad
\xi_N=-\sum_{n=1}^N\tpsi_{-2(n-1)}\tpsi_{-2(n-1)}^\ast+
\sum_{n=1}^N\tpsi_{2n}^\ast\tpsi_{2n}.
\non
\ena
They act on the direct sum of the spaces
\bea
&&
\tH_{2m}^{(N)}=
\sum_{-N+1<r_1<\cdots<r_k\leq 0<s_1<\cdots<s_{k+m}\leq N}
\mathbb{C}
\tpsi_{2r_1}\cdots \tpsi_{2r_k}
\tpsi^\ast_{2s_1}\cdots \tpsi^\ast_{2s_{k+m}}|0>,
\non
\ena
and satisfy the relations
\bea
&&
[\eta_N,\omega_N]=\xi_N,
\quad
[\xi_N,\eta_N]=2\eta_N,
\quad
[\xi_N,\omega_N]=-2\omega_N.
\non
\ena
Thus $\omega_N$, $\eta_N$ and $\xi_N$ determine the action of $sl_2$ on 
$\oplus_{m\in \mathbb{Z}}\tH_{2m}^{(N)}$ by the correpondence
\bea
&&
e=\eta_N,
\quad
f=\omega_N,
\quad
h=\xi_N.
\non
\ena
We have $h=m$ on $\tH_{2m}^{(N)}$. Thus by the representation theory of
$sl_2$, $\omega_N$ is injective at $\tH_{2m}^{(N)}$, $m\geq1$. Since
$v\in \tH_{2m}^{(N)}$, $v=0$.
\qed

Let us write $D=q^d\in \uq$. Then $d$ can be considered as a degree operator
of $V(\Lambda_i)$. The grading of $V(\Lambda_i)$ induces those of
$\vis1$ and $H_{2m}$.
Since $x^{-}_0$ and $(x^{-}_0)^{(2)}$ are homogeneous,
the quotient space $A_{2m}$ is also graded. 
In general, for a graded vector space $V$ by the degree operator $d$ such that
each homogeneous subspace is finite dimensional, we define its character by
\bea
&&
\ch\, V=\tr_{V}(\,p^{-d}\,).
\non
\ena
We remark that, with this definition of the degree, we have
\bea
&&
\deg\, \psi_{2n}=2n-1,
\quad
\deg\, \psi_{2n}^\ast=-(2n-1),
\quad
\deg\, \,|2m>=-m^2.
\non
\ena
Then
\bea
&&
\deg\, \tpsi_{2n}=2n-1,
\quad
\deg\, \tpsi_{2n}^\ast=-(2n-1),
\non
\ena
and
\bea
&&
\ch\, \tH_{2m}=\ch\, H_{2m}.
\label{h=th}
\ena

\begin{cor}
\bea
&&
\ch\, A_{2m}=
\frac{p^{m^2}(1-p^{2m+1})}{(p:p)_\infty}.
\non
\ena
\end{cor}

\noindent
{\it Proof.} Since
\bea
&&
H_{2m}\simeq \mathbb{C}[b_{-2},b_{-4},\ldots]\e^{\Lambda_0+m\alpha_1},
\non
\ena
we have, by (\ref{h=th}),
\bea
&&
\ch\, \tH_{2m}=\frac{p^{m^2}}{(p^2:p^2)_\infty}.
\non
\ena
By Lemma \ref{lem2} $\omega$ is injective. Thus
\bea
&&
\ch\, W_{2m}=\ch\, \tH_{2m}-\ch\, \tH_{2m+4},
\non
\ena
and
\bea
&&
\ch\, A_{2m}=
\sum_{n=m}^\infty (-1)^{n-m}
\ch\Big(
{\cal D}\otimes W_{2m}
\Big)
=\ch\, {\cal D}(\ch\, \tH_{2m}-\ch\, \tH_{2m+2}).
\non
\ena
This proves the corollary.
\qed
\vskip2mm

\noindent
{\bf Remark.} The character of $A_{2m}$ is calculated 
in \cite{N2,FJKMMT} in a so called fermionic form. 
Combining the identity of $q$-series \cite{Me} with it the character formula in the corollary is proved. Here we have shown that the character formula 
in this bosonic form can be derived from the resolution of 
$A_{2m}$ independently to other results.
\vskip2mm

\begin{cor}\label{highest}
We have the isomorphism
\bea
&&
\frac{A_{2m}}{\sum_{s=1}^\infty b_{-(2s-1)}A_{2m}}
\simeq W_{2m}.
\non
\ena
\end{cor}

\section{De Rham type cohomologies}
Let us set
\bea
&&
v_\phi=dt_{-1}\wedge dt_{-3}\wedge dt_{-5}\wedge\cdots.
\non
\ena
For $0>i_1>\cdots>i_l$, $i_k$'s being odd,
we define $v_{i_1,\ldots,i_l}$ to be the formal infinite wedges obtained from 
$v_\phi$ by removing $dt_{i_1}$,...,$dt_{i_l}$.
We define $\Omega^{\frac{\infty}{2}-p}$, $p\geq 0$ by
\bea
&&
\Omega^{\frac{\infty}{2}}=\mathbb{C}v_\phi,
\non
\\
&&
\Omega^{\frac{\infty}{2}-p}=\sum_{0>i_1>\cdots>i_p, \text{$i_k$ odd}}
\mathbb{C}v_{i_1,\ldots,i_p}
\quad
\text{for $p\geq 1$}
\non
\ena
and 
set
\bea
&&
C_{2m}^{\frac{\infty}{2}-p}=A_{2m}\otimes \Omega^{\frac{\infty}{2}-p}.
\non
\ena
Then by defining
\bea
&&
d=\sum_{s=1}^\infty b_{-(2s-1)}\otimes dt_{-(2s-1)},
\non
\ena
to be a differential $(C_{2m}^{\frac{\infty}{2}-\cdot},d)$ becomes a complex.
We denote its cohomology group at $C_{2m}^{\frac{\infty}{2}-p}$ by 
$H_{2m}^{\frac{\infty}{2}-p}$.
Notice that 
\bea
&&
H_{2m}^{\frac{\infty}{2}}\simeq 
\frac{A_{2m}}{\sum_{s=1}^\infty b_{-(2s-1)}A_{2m}},
\non
\ena
which is described in Corollary \ref{highest}.
Other cohomology groups are described by

\begin{prop}\label{de Rham}
\bea
&&
H_{2m}^{\frac{\infty}{2}-p}
\simeq 
W_{2(m+p)}
\quad\quad
\text{for} \quad p\geq 0.
\non
\ena
\end{prop}

As a consequence of this proposition we have
\bea
&&
\ch\, H_{2m}^{\frac{\infty}{2}-p}
=\frac{p^{(m+p)^2}(1-p^{4(m+p+1)})}{(p^2:p^2)_\infty}.
\non
\ena

\vskip2mm
\noindent
{\it Proof of Proposition \ref{de Rham}}
\par
Consider two chain complexes of ${\cal D}$-modules
\bea
&&
K_p={\cal D}\otimes \Omega^{\frac{\infty}{2}-p},
\quad
d^K=\sum_{n=1}^\infty b_{-(2n-1)}\otimes dt_{-(2n-1)},
\non
\\
&&
L_p={\cal D}\otimes W_{2(m+p)},
\quad
d^L=x^{-}_0.
\non
\ena
We set $K_p=L_p=0$ for $p<0$.
Let 
\bea
&&
M_{p,q}=K_p\otimes_{{\cal D}} L_q
\non
\ena
be the double complex obtained from $K$ and $L$. Here $p$ denotes the row index and $q$ denotes the column index. 
We denote by $M=(M_n)$ the total complex of $(M_{p,q})$.
The bordered chain complexes of $(M_{p,q})$
in the vertical and horizontal directions are
\bea
&&
H_0(M^{p,\cdot})\simeq C_{2m}^{\frac{\infty}{2}-p},
\label{border1}
\\
&&
H_0(M^{\cdot,q})\simeq W_{2(m+q)},
\label{border2}
\ena
respectively, where the symbol $H_0$ denotes the $0$-th homology group 
of a comolex.  Notice that all maps in the complex (\ref{border2}) are zero.
Therefore $q$-th homology group of (\ref{border2}) are given by
\bea
&&
H_q(W_{2(m+\cdot)})=W_{2(m+q)}.
\non
\ena
By the definition $p$-th homology group of the complex (\ref{border1}) is
\bea
&&
H_p(C^{\frac{\infty}{2}-\cdot})=H^{\frac{\infty}{2}-p}_{2m}.
\non
\ena
Since two complexes 
\bea
&&
\cdots \lar M_{1,q} \lar M_{0,q} \lar W_{2(m+q)} \lar 0,
\non
\\
&&
\cdots \lar M_{p,1} \lar M_{p,0} \lar C_{2m}^{\frac{\infty}{2}-p} \lar 0,
\non
\ena
are exact for $p,q\geq 0$, we have
\bea
&&
H_n(W_{2(m+\cdot)})\simeq H_n(M) \simeq H_n(C_{2m}^{\frac{\infty}{2}-\cdot}).
\non
\ena
\qed

\section{The structure of $A_{2m+1}$}
\par
The strategy is similar to the case of $A_{2m}$, that is, we make a base change of fermions and describe $A_{2m+1}$ using the Fock space of 
new fermions. The structure of $A_{2m+1}$ is simpler than that of 
$A_{2m}$.

The component forms of $x^{-}_0$ and $(x^{-}_0)^{(2)}$ are given by
\bea
&&
x^{-}_0=\sum_{n=0}^\infty T_{-2n} \psi_{2n+1},
\non
\\
&&
(x^{-}_0)^{(2)}=
\frac{-i}{2}
\sum_{n=1}^\infty
\Big(
-2\psi_{-2n+1}+\sum_{l=1}^\infty Q_{n-1,l}(T)\psi_{2(-n+l)+1}
\Big)
\psi_{2n+1},
\non
\ena
where
\bea
&&
Q_{n,l}(T)=
\sum_{n_1+n_2=l,0\leq n_1,0\leq n_2\leq n}
(-2T_{-2n_1}T_{-2n_2}+T_{-(2n_1-1)}T_{-(2n_2+1)}).
\non
\ena
We set
\bea
&&
\tpsi_{2n+1}=\psi_{2n+1}\quad \text{for $n\geq 0$},
\non
\\
&&
\tpsi_{-(2n+1)}=-2\psi_{-(2n+1)}+\sum_{l=1}^\infty Q_{n,l}(T)\psi_{2(-n+l)-1}
\quad \text{for $n\geq 0$}.
\non
\ena
Then we have
\bea
&&
x^{-}_0=\sum_{n=0}^\infty T_{-2n} \tpsi_{2n+1},
\non
\\
&&
(x^{-}_0)^{(2)}=
\frac{-i}{2}
\sum_{n=1}^\infty\tpsi_{-(2n-1)}\tpsi_{2n+1}.
\non
\ena
We define $\tpsi_{2n+1}^\ast$, $n\in \mathbb{Z}$ similarly to the case of 
$A_{2m}$. Again the vacuums are invariant by this change. Thus the Fock
space $\tH_{2m+1}$ is defined over the same vacuum as fermions 
$\psi_{2n+1}$, $\psi^\ast_{2n+1}$. 

Let us define the subspace of $\tH_{2m+1}$ by
\bea
&&
\tH_{2m+1,0}=
\sum 
\mathbb{C}
\tpsi_{2r_1+1}\cdots\tpsi_{2r_k+1}
\tpsi^\ast_{2s_1+1}\cdots\tpsi_{2s_{k+m}+1}^\ast|1>,
\non
\ena
where
$r_1<\cdots<r_k<0<s_1<\cdots<s_{k+m}$. 
Notice that $\tpsi_1$ is absent in $\tH_{2m+1,0}$.
Since $(x^{-}_0)^{(2)}$ does not contain $\tpsi_1$,
\bea
&&
(x^{-}_0)^{(2)}\tH_{2m+1,0}\subset \tH_{2m-3,0}.
\non
\ena
We set
\bea
&&
W_{2n+1,0}=\frac{\tH_{2n+1,0}}{(x^{-}_0)^{(2)}\tH_{2n+5,0}}.
\non
\ena
The degree induced from $V_{\sqrt{-1}}(\Lambda_1)$ is given by
\bea
&&
\deg\, \tpsi_{2n+1}=2n,
\quad
\deg\, \tpsi_{2n+1}^\ast=-2n,
\quad
\deg\,|2m+1>=-m^2-m.
\non
\ena
\begin{theorem}\label{oddspace}
(i) The space $A_{2m+1}$ becomes a free ${\cal D}$-module as
\bea
&&
A_{2m+1}\simeq
{\cal D}\otimes W_{2m+1,0}.
\non
\ena

\noindent
(ii)
\bea
&&
\ch\, A_{2m+1}=\frac{p^{m(m+1)}(1-p^{2m+2})}{(p:p)_\infty}.
\non
\ena
\end{theorem}

\begin{cor}\label{oddcohomology}
\bea
&&
\frac{A_{2m+1}}{\sum_{s=1}^\infty b_{-(2s-1)}A_{2m+1}}
\simeq
W_{2m+1,0}.
\non
\ena
\end{cor}

\vskip2mm
\noindent
{\it Proof of Theorem \ref{oddspace}}.
\par
Since
\bea
&&
x^{-}_0=\tpsi_1+T_2\tpsi_3+T_4\tpsi_5+\cdots,
\non
\ena
we have the isomorphism of ${\cal D}$-modules
\bea
&&
{\cal D}\otimes \tH_{2m+1,0}
\simeq
\frac{{\cal D}\otimes \tH_{2m+1}}{x^{-}_0({\cal D}\otimes \tH_{2m+1})}.
\non
\ena
Thus the natural map
\bea
&&
{\cal D}\otimes \tH_{2m+1,0}
\lar
A^{(0)}_{2m+1}
\non
\ena
is surjective and its kernel is given by
\bea
&&
\frac{x^{-}_0({\cal D}\otimes \tH_{2m+1})+
(x^{-}_0)^{(2)}({\cal D}\otimes \tH_{2m+5})}
{x^{-}_0({\cal D}\otimes \tH_{2m+1})}
\simeq 
(x^{-}_0)^{(2)}({\cal D}\otimes \tH_{2m+5,0}),
\non
\ena
which proves (i).

Next let us prove (ii).
Set
\bea
&&
\omega_N=-\sum_{n=1}^N \tpsi_{-(2n-1)}\tpsi_{2n+1},
\non
\\
&&
\eta_N=\sum_{n=1}^N \tpsi^\ast_{-(2n-1)}\tpsi^\ast_{2n+1},
\non
\\
&&
\xi_N=-\sum_{n=1}^N\tpsi_{-(2n-1)}\tpsi^\ast_{-(2n-1)}
+\sum_{n=1}^N\tpsi^\ast_{2n+1}\tpsi_{2n+1}.
\non
\ena
Then they satisfy the relations of $sl_2$ as in the case of $A_{2m}$.
It follows that $(x^{-}_0)^{(2)}$ is injective on $\tH_{2m+1,0}$, $m\geq 1$.
Using the arguments of \cite{BBS} we have
\bea
\ch\, W_{2m+1,0}
&=&\ch\, \tH_{2m+1,0}-\ch\, \tH_{2m+5,0}
\non
\\
&=&
(1-p^{2m+2})\ch\, \tH_{2m+1}.
\non
\ena
Thus (ii) follows from
\bea
&&
\ch\, \tH_{2m+1}=\frac{p^{m(m+1)}}{(p^2:p^2)_\infty}.
\non
\ena
\qed

We define de Rham type complex $(C_{2m+1}^{\frac{\infty}{2}-p},d)$ 
and its cohomology $H^{\frac{\infty}{2}-p}_{2m+1}$for $A_{2m+1}\,$ in a similar manner to the case of $A_{2m}$. Then

\begin{prop}
$$
H^{\frac{\infty}{2}-p}_{2m+1}\simeq
\left\{
\begin{array}{rl}
W_{2m+1,0},& \quad \text{for $p=0$}\\
0,&\quad \text{for $p\geq 1$}.
\end{array}\right.
$$
\end{prop}

The proof of this proposition is similar to that of Proposition
\ref{de Rham} and we leave it to the reader.

\section{Action of symplectic group on cohomologies}

Let $Sp(2n)\subset GL(2n,\mathbb{C})$ be the symplectic group.
For $n<n'$ we have the embedding
\bea
&&
GL(2n,\mathbb{C})\subset GL(2n',\mathbb{C}),
\non
\\
&&
A\mapsto
\left[
\begin{array}{ccc}
1_{n'-n}&\quad&\quad\\
\quad&A&\quad\\
\quad&\quad&1_{n'-n}
\end{array}
\right],
\non
\ena
where $1_{r}$ denotes the $r$ by $r$ unit matrix.
We define the group $Sp(2\infty)$ as the inductive limit of $Sp(2n)$
with respect to this embedding. In this section we shall show that $W_{2m}$
and $W_{2m+1,0}$ are irreducible representations of $Sp(2\infty)$.

Let us consider $W_{2m}$ first. Set 
\bea
&&
\alpha_n=\tpsi^\ast_{-2(n-1)},
\quad
\beta_n=\tpsi^\ast_{2n},
\quad
\alpha_n^\ast=\tpsi_{-2(n-1)},
\quad
\beta_n^\ast=\tpsi_{2n},
\quad
n\geq 1.
\non
\ena
Then
\bea
&&
\omega=-\sum_{n=1}^\infty \alpha_n^\ast \beta_n^\ast,
\quad
\eta=\sum_{n=1}^\infty \alpha_n\beta_n.
\label{eq7-1}
\ena
Let
\bea
&&
V_N=\oplus_{i=1}^N\mathbb{C}\alpha_i\oplus_{i=1}^N\mathbb{C}\beta_i.
\label{eq70}
\ena
We consider $V_N$ as the vector representation of $Sp(2N)$.
Then for $k\geq 1$
\bea
&&
M^{(N)}_k:=\frac{\wedge^k V_N}{\eta_N \wedge^{k-2}V_N},
\label{eq71}
\ena
is isomorphic to the $k$-th fundamental representation of $Sp(2N)$,
where 
\bea
&&
\eta_N=\sum_{i=1}^N \alpha_i \wedge \beta_i.
\non
\ena
To each $m\geq 0$ we associate $M^{(N)}_{N-m}$, $N\geq m$.
For $N<N'$ we define a map
\bea
&&
\wedge^{N-m} V_N \lar \wedge^{N'-m} V_{N'},
\non
\\
&&
v \mapsto v \wedge \alpha_{N+1} \wedge \cdots \wedge \alpha_{N'}.
\label{eq72}
\ena
It induces a map
\bea
&&
M^{(N)}_{N-m} \lar M^{(N')}_{N'-m}.
\label{eq73}
\ena

\begin{lemma}
The map (\ref{eq73}) is injective.
\end{lemma}

\vskip2mm
\noindent
{\it Proof.}
It is sufficient to prove the lemma for $N'=N+1$.
Suppose that $v\in \wedge^{N-m} V_N$ satisfy
\bea
&&
v\wedge \alpha_{N+1} =\eta_{N+1}\wedge w,
\quad
w\in \wedge^{N-m-1} V_{N+1}.
\label{eq74}
\ena
Let us write
\bea
&&
w=w'+w''\wedge \alpha_{N+1},
\quad
w'\in \wedge^{N-m-1} V_{N+1},
\quad
w''\in \wedge^{N-m-2} V_{N+1},
\non
\ena
where $w'$ and $w''$ do not contain $\alpha_{N+1}$.
Then
\bea
&&
\eta_{N+1} \wedge w=\eta_N\wedge w'+
(-w' \wedge \beta_{N+1}+\eta_N \wedge w'') \wedge \alpha_{N+1}.
\non
\ena
By (\ref{eq74}) we have
\bea
&&
\eta_N\wedge w'=0.
\non
\ena
By the representation theory of $sl_2$ $\eta_N$ is injective on 
$\wedge^{N-m-1} V_{N}\oplus \beta_{N+1}\wedge^{N-m-2} V_{N}$. 
Thus $w'=0$ and we have
\bea
&&
v \wedge \alpha_{N+1}=\eta_{N} \wedge w'' \wedge \alpha_{N+1}.
\non
\ena
This shows $v=\eta_n \wedge w''$.
\qed
\vskip2mm

We denote $M_{\infty-m}$ the inductive limit of $M^{(N)}_{N-m}$,
\bea
&&
M_{\infty-m}=\lim_{\lar}M^{(N)}_{N-m}.
\non
\ena
For $N<N'$ the subgroup $Sp(2N)$ in $Sp(2N')$ fixes 
$\alpha_{N+1}$,...,$\alpha_{N'}$, $\beta_{N+1}$,...,$\beta_{N'}$.
Therefore $Sp(2\infty)$ acts on $M_{\infty-m}$. It is straightforward to
check that this representation is irreducible.

\begin{prop}
$$W_{2m}\simeq M_{\infty-m}.$$
\end{prop}
\vskip2mm
\noindent
{\it Proof.}
By the representation theory of $sl_2$ we have the isomorphisms,
\bea
&&
\omega^m:\, \tH_{2m}\simeq \tH_{-2m},
\non
\\
&&
\omega^{m+1} \tH_{2m+4}\simeq \eta \tH_{-2m-4}.
\non
\ena
Thus
\bea
&&
W_{2m}=
\frac{\tH_{2m}}{\omega \tH_{2m+4}}
\simeq
\frac{\tH_{-2m}}{\eta \tH_{-2m-4}}.
\non
\ena
We shall show that $W_{2m}$ is an inductive limit of the subspaces
isomorphic to $M^{(N)}_{N-m}$.

We set, for $N\geq m$,
\bea
&&
\tH_{-2m}(N)=\sum_{k=0}^{N-m}\sum_{}\mathbb{C}
\beta_{i_1}\cdots\beta_{i_k}\alpha_{j_1}\cdots \alpha_{j_{N-m-k}}|-2N>,
\non
\ena
where the second summation is taken for all
\bea
&&
1\leq i_1<\cdots<i_k\leq N,
\quad
1\leq j_1<\cdots<j_{N-m-k}\leq N.
\non
\ena
For $N<N'$ we have the inclusion
\bea
&&
\tH_{-2m}(N)\subset \tH_{-2m}(N'),
\label{eq75}
\\
&&
x|-2N>= x\alpha_{N+1}\cdots\alpha_{N'}|-2N'>.
\non
\ena
Thus $\{\tH_{-2m}(N)\}$ defines an increasing filtration and satisfy
\bea
&&
\tH_{-2m}=\cup_{N=m}^\infty \tH_{-2m}(N).
\non
\ena
There is an isomorphism,
\bea
&&
\tH_{-2m}(N) \simeq \wedge^{N-m} V_N,
\label{eq76}
\\
&&
\beta_{i_1}\cdots\beta_{i_k}\alpha_{j_1}\cdots \alpha_{j_{N-m-k}}|-2N>
\mapsto
\beta_{i_1}\wedge \cdots \wedge \alpha_{j_{N-m-k}}.
\non
\ena
It induces the isomorphism
\bea
&&
\frac{\tH_{-2m}(N)}{\eta_N \tH_{-2m-4}(N)}
\simeq
\frac{\wedge^{N-m} V_N}{\eta_N \wedge^{N-m-2} V_N}
=
M^{(N)}_{N-m}.
\non
\ena
By the isomorphism (\ref{eq76}) the inclusion (\ref{eq75}) is transformed to
 the map (\ref{eq72}). Thus we have
\bea
&&
\frac{\tH_{-2m}}{\eta \tH_{-2m-4}}
=
\lim_{\lar}\frac{\tH_{-2m}(N)}{\eta_N \tH_{-2m-4}(N)}
\simeq 
M_{\infty-m}.
\non
\ena
\qed

As to $W_{2m+1,0}$ we have

\begin{prop}
$$W_{2m+1,0}\simeq M_{\infty-m-1}.$$
\end{prop}
\vskip2mm
\noindent
{\it Proof.}
In this case we set
\bea
&&
\alpha_n=\tpsi^\ast_{-(2n-1)},
\quad
\beta_n=\tpsi^\ast_{2n+1},
\quad
\alpha_n^\ast=\tpsi_{-(2n-1)},
\quad
\beta_n^\ast=\tpsi_{2n+1},
\quad
n\geq 1.
\non
\ena
Then $\omega$ and $\eta$ are given by (\ref{eq7-1}).
Define $V_N$ and $M^{(N)}_k$ by (\ref{eq70}) and (\ref{eq71}) respectively.
As in the case of $\tH_{2m}$ we have the isomorphisms
\bea
&&
\omega^{m+1}:\, \tH_{2m+1,0}\simeq \tH_{-2m-1,0},
\non
\\
&&
\omega^{m+2}\tH_{2m+5,0}\simeq \eta \tH_{-2m-5,0},
\non
\ena
and consequently
\bea
&&
W_{2m+1,0}\simeq \frac{\tH_{-2m-1,0}}{\eta \tH_{-2m-5,0}}.
\non
\ena
By the definition
\bea
&&
\tH_{-2m-1,0}=
\sum_{k=m+1}^\infty 
\sum_{}
\mathbb{C}
\beta_{i_1}\cdots\beta_{i_{k-m-1}}
\alpha_{j_1}^\ast\cdots \alpha_{j_k}^\ast|1>,
\non
\ena
where the second summation is taken for all
\bea
&&
1\leq i_1<\cdots<i_{k-m-1},
\quad
1\leq j_1<\cdots<j_k.
\non
\ena
We identify $\tH_{-2m-1,0}$ as a subspace of $\tH_{-2m-3}$ by
\bea
&&
x|1>=x\tpsi_1^\ast|-1>\mapsto x|-1>,
\label{eq77}
\ena
that is, by removing $\tpsi_1^\ast$. Since $x$ does not contain $\tpsi_1$, this map is well-defined. Moreover the map (\ref{eq77}) preserves the degree because
$\deg\, \tpsi_1^\ast=0$. We set
\bea
&&
\tH_{-2m-1,0}(N)=
\sum_{k=0}^{N-m-1}
\sum_{}
\mathbb{C}
\beta_{i_1}\cdots\beta_{i_{k}}
\alpha_{j_1}\cdots \alpha_{j_{N-m-1-k}}|-(2N+1)>,
\non
\ena
where the second summation is taken for all
\bea
&&
1\leq i_1<\cdots<i_{k}\leq N,
\quad
1\leq j_1<\cdots<j_{N-m-1-k}\leq N.
\non
\ena
Then 
\bea
&&
\tH_{-2m-1,0}(N)\simeq \wedge^{N-m-1} V_N.
\label{eq78}
\ena
For $N<N'$ we have the natural inclusion map
\bea
&&
\tH_{-2m-1,0}(N)\lar \tH_{-2m-1,0}(N'),
\non
\\
&&
x|-(2N+1)>\mapsto x\alpha_{N+1}\cdots \alpha_{N'}|-(2N'+1)>.
\non
\ena
This map is transplanted to (\ref{eq72}) by the isomorphism (\ref{eq78}).
Thus we have
\bea
&&
\frac{\tH_{-2m-1,0}}{\eta \tH_{-2m-5,0}}
\simeq
\lim_{\lar}
\frac{\tH_{-2m-1,0}(N)}{\eta_N \tH_{-2m-5,0}(N)}
\simeq 
M_{\infty-m-1}.
\non
\ena
\qed

\end{document}